\documentclass[9pt,twoside,reqno,draft]{amsart}

\newcommand{\vs}{\vspace}
\newcommand{\bc}{\begin{center}}
\newcommand{\ec}{\end{center}}
\newcommand{\bt}{\begin{theorem}}
\newcommand{\et}{\end{theorem}}
\newcommand{\bl}{\begin{lemma}}
\newcommand{\el}{\end{lemma}}
\newcommand{\br}{\begin{remark}}
\newcommand{\er}{\end{remark}}
\newcommand{\bco}{\begin{corollary}}
\newcommand{\eco}{\end{corollary}}
\newcommand{\bp}{\begin{proof}}
\newcommand{\ep}{\end{proof}}
\newcommand{\be}{\begin{example}}
\newcommand{\ee}{\end{example}}
\newcommand{\bd}{\begin{definition}}
\newcommand{\ed}{\end{definition}}
\newcommand{\ba}{\begin{align}}
\newcommand{\ea}{\end{align}}
\newcommand{\bpr}{\begin{proposition}}
\newcommand{\epr}{\end{proposition}}
 
\usepackage{amsthm,amssymb,amsmath,color,amsfonts,fancybox,fancyhdr,mathrsfs,fmultico}

\newtheorem{theorem}{Theorem}[section]
\newtheorem{definition}{Definition}[section]
\newtheorem{example}{Example}[section]
\newtheorem{lemma}{Lemma}[section]
\newtheorem{remark}{Remark}[section]
\newtheorem{corollary}{Corollary}[section]

\begin{document}
\setcounter{page}{1}

\leftline{\footnotesize {\bf\em My Paper ,2012}}
\leftline{\footnotesize
masadi@azu.ac.ir}

\vs*{2.7cm}

\title[An Inequality ...]{\large An Inequality of Uniformly Continuous Functions in Normed Spaces}
\author[Mehdi Asadi]{\sc Mehdi Asadi$^{1}$}
\date{}
\maketitle

\vs*{-0.5cm}

\bc
{\footnotesize $^{}\footnote{Corresponding author.\quad Fax:+98-241-4220030.}$  Department of Mathematics, Zanjan Branch, Islamic Azad University,  Zanjan, Iran\\
masadi@azu.ac.ir\\
\medskip
}
\ec


{\footnotesize
\noindent
{\bf Abstract.}
We obtain an interesting inequalities for uniformly continuous functions in the normed spaces: $\|f(x)\|\leq a\|x\|+b$ for some $a,b> 0$.\newline

\noindent
\bigskip
\section{Introduction}
We start first with an inequality on normed spaces by
uniformly continuous functions and then this shall obtain
inequality for uniformly continuous functions on $\mathbb{R}$ as simple case.
\section{Main result}
\begin{theorem}
If $X,Y$ be two linear normed spaces on $\mathbb{R}$ and
$f$ uniformly continuous function on $X$, then there exist two real numbers $a,b>0$ such that $\|f(x)\|\leq a\|x\|+b$ for all $x\in X.$
\end{theorem}
{\bf Proof.} Since $f$ is uniformly continuous function on $X$, so
$$\exists \delta>0\quad \forall x,y\in X\quad \|x-y\|<\delta\Rightarrow \|f(x)-f(y)\|<1.$$
Let $x\in X$ be arbitrary and fixed, put
$0<d_0:=\frac{\delta}{2}<\delta$ choose $0\not =x_0\in X$ such
that $0<\|x_0\|<d_0$, now $0<d_0-\|x_0\|$, choose $x_0\not =x_1\in
X$ and $ x_1\not= 0$ such that $0<\|x_1\|<d_0-\|x_0\|$, therefore
$0<\|x_1-x_0\|\leq\|x_1\|+\|x_0\|<d_0,$ so $ \|x_1-x_0\|<\delta.$
Similar to since $0<\|x_1\|<d_0-\|x_0\|<d_0$, this implies that
$0<d_0-\|x_1\|$, then choose $x_1\not =x_2\in X$ and $ x_2\not= 0$ such
that $0<\|x_2\|<d_0-\|x_1\|$, therefore
$$0<\|x_2-x_1\|\leq\|x_2\|+\|x_1\|<d_0<\delta,$$
by continuing this procedure up to where $N$ be a smallest natural
number such that $N-1<\frac{2}{\delta}\|x\|\leq N$, so there
exists $0\not=x_N\in X$ and $ x_N\not =x_{N-1}$ such that
$0<\|x_N\|<d_0-\|x_{N-1}\|$ and
$$
0<\|x_N-x_{N-1}\|<\delta,\quad \|x_N\|<d_0,$$
and also
$$\|x-x_N\|\leq \|x\|+\|x_N\|\leq\frac{N\delta}{2}+d_0=\frac{(N+1)\delta}{2}.$$
Now put $y_0:=x\in X$ and take
$y_1:=\frac{N}{N+1}x+\frac{1}{N+1}x_N$ so $y_1\in X$ and
$$\|y_0-y_1\|=\frac{\|x-x_n\|}{N+1}<\frac{\delta}{2},$$
and take $y_2:=2y_1-y_0$ then
$$\|y_2-y_1\|=\|y_1-y_0\|=\|y_1-x\|<\frac{\delta}{2},$$
so by this way  $y_N:=2y_{N-1}-y_{N-2}$ and $y_N:=x_N$, therefore
$$\|x_N-y_{N-1}\|=\|y_N-y_{N-1}\|=\|y_{N-1}-y_{N-2}\|<\frac{\delta}{2},$$
thus $y_i\in X, (0\leq i\leq N)$ and
$\|y_i-y_{i-1}\|<\frac{\delta}{2}$ so $$\|f(y_i)-f(y_{i-1})\|<1,$$
and since $\|x_1\|<\delta$ and $\|x_i-x_{i+1}\|<\delta$ for $1\leq i\leq N-1,$
hence we have
\begin{eqnarray}
\|f(0)-f(x)\|&\leq&\|f(0)-f(x_1)\|+\sum_{i=1}^{N-1}\|f(x_i)-f(x_{i+1})\|+\|f(x_N)-f(x)\|,\nonumber\\
&\leq& 1+1+\cdots+1+\|f(x_N)-f(x)\|=N+\|f(x_N)-f(x)\|,\nonumber\\
&\leq &N+\|f(x)-f(y_1)\|+\|f(y_1)-f(y_2)\|+\cdots+\|f(y_{N-1})-f(x_N)\|,\nonumber\\
&=&N+\|f(y_0)-f(y_1)\|+\|f(y_1)-f(y_2)\|+\cdots+\|f(y_{N-1})-f(y_N)\|,\nonumber\\
&<& N+N\times 1=2N.\nonumber
\end{eqnarray}
Now from $N-1<\frac{2}{\delta}\|x\|\leq N$ we obtain
$N<1+\frac{2}{\delta}\|x\|$ consequentially
\begin{eqnarray}
\|f(x)\|&\leq&\|f(x)-f(0)\|+\|f(0)\|,\nonumber\\
&\leq&2N+\|f(0)\|,\nonumber\\
&\leq &\frac{4}{\delta}\|x\|+2+\|f(0)\|=a\|x\|+b,\nonumber
\end{eqnarray}
where $a:=\frac{4}{\delta}, b:=2+\|f(0)\|.$$\hfill\qed$

\begin{corollary}(See ~\cite{cite1}, page 47)
If $f:\mathbb{R} \rightarrow \mathbb{R}$ be uniformly continuous function,
then there exist $a,b>0$ such that
$|f(x)|\leq a|x|+b,$
for all $x\in X.$
\end{corollary}


\end{document}